\documentclass{amsproc}
\usepackage{aspmproc}
\usepackage[mathscr]{eucal}
\usepackage{amssymb}
\def\Bbb{\mathbb}
\def\frak{\mathfrak}

\newenvironment{pf*}[1]{\proof[#1]}{\endproof}
\newcommand{\rom}{\textup}
\newtheorem{MainTheorem}{Theorem}
\newtheorem{Corollary}[equation]{Corollary}
\newtheorem{Lemma}[equation]{Lemma}
\newtheorem{Proposition}[equation]{Proposition}

\theoremstyle{definition}
\newtheorem{Definition}[equation]{Definition}

\theoremstyle{remark}
\newtheorem{Remark}[equation]{Remark}

\numberwithin{equation}{section}

\newcommand{\thmref}[1]{Theorem~\ref{#1}}

\newcommand{\lemref}[1]{Lemma~\ref{#1}}
\newcommand{\propref}[1]{Proposition~\ref{#1}}
\newcommand{\corref}[1]{Corollary~\ref{#1}}
\newcommand{\subsecref}[1]{\S\ref{#1}}

\newcommand{\remref}[1]{Remark~\ref{#1}}
\newcommand{\defeq}{\overset{\operatorname{\scriptstyle def.}}{=}}
\newcommand{\C}{{\Bbb C}}
\newcommand{\Z}{{\Bbb Z}}
\newcommand{\Q}{{\Bbb Q}}

\newcommand{\GL}{\operatorname{GL}}

\newcommand{\g}{{\frak g}}

\newcommand{\Hom}{\operatorname{Hom}}

\newcommand{\Ker}{\operatorname{Ker}}

\newcommand{\ve}{\varepsilon}
\newcommand{\Uq}{{\mathbf U}_q} 

\newcommand{\Um}{\Tilde{\mathbf U}_q(\ag)} 
\newcommand{\Ul}{{\mathbf U}_q({\mathbf L}{\mathfrak g})} 

\newcommand{\Uli}{{\mathbf U}^{\Z}_q({\mathbf L}\mathfrak g)}

\newcommand{\ag}{\widehat{\g}}
\newcommand{\wt}{\operatorname{wt}}
\newcommand{\bc}{\mathbf c}
\newcommand{\te}{{\widetilde e}}
\newcommand{\tf}{{\widetilde f}}
\newcommand{\B}{\mathcal B} 
\newcommand{\La}{{\mathscr L}} 
\newcommand{\aW}{\widehat W} 
\newcommand{\eW}{\widetilde W} 
\newcommand{\Ua}{{\mathbf U}_q(\ag)} 
\newcommand{\Ui}{{\mathbf U}^{\Z}_q(\widehat{\mathfrak g})}
\newcommand{\Uqm}{\widetilde{\mathbf U}_q(\ag)} 
\newcommand{\A}{\mathbf A}
\newcommand{\aR}{\mathscr R}

\makeatletter
\xdef\widecheck#1{\noexpand\@mathmeasure\z@\textstyle{#1}%
  \noexpand\ifdim\noexpand\wdz@>\tw@ em%
  \mathaccent"0\hexnumber@\symAMSb 5B{#1}%
  \noexpand\else\mathaccent"0371{#1}\noexpand\fi}
\makeatother

\title{Extremal weight modules of quantum affine algebras}
\author[H. Nakajima]{Hiraku Nakajima}
\address{Department of Mathematics\\
Kyoto University\\
Kyoto 606-8502\\
Japan}

\begin{document}
\begin{abstract}
  Let $\ag$ be an affine Lie algebra, and let $\Ua$ be the quantum
  affine algebra introduced by Drinfeld and Jimbo. In \cite{Kas94}
  Kashiwara introduced a $\Ua$-module $V(\lambda)$, having a global
  crystal base for an integrable weight $\lambda$ of level $0$. We
  call it an {\it extremal weight module}.
  It is isomorphic to the Weyl module introduced by Chari-Pressley
  \cite{CP00}.
  In \cite[\S13]{Kas00} Kashiwara gave a conjecture on the structure
  of extremal weight modules. We prove his conjecture when $\ag$ is
  an untwisted affine Lie algebra of a simple Lie algebra $\g$ of
  type $ADE$, using a result of Beck-Chari-Pressley \cite{BCP}. As a
  by-product, we also show that the extremal weight module is
  isomorphic to a universal standard module, defined via quiver
  varieties by the author \cite{Nak00,Nak01}. This result was
  conjectured by Varagnolo-Vasserot \cite{VV} and Chari-Pressley
  \cite{CP00} in a less precise form.
  Furthermore, we give a characterization of global crystal bases by
  an almost orthogonality propery, as in the case of global crystal
  base of highest weight modules.
\end{abstract}
\maketitle

\section{Introduction}

In the conference, I gave a survey on quiver varieties and finite
dimensional representations of quantum affine algebras. Since I
already wrote a survey article \cite{Na-qchar} on this subject, I will
discuss a different one in this paper. But it is related to my talks
since I will study extremal weight modules which turn out to be
isomorphic to universal standard modules, which was one of the main
objects in my talk.

Let us describe Kashiwara's conjecture \cite[\S13]{Kas00} on extremal
weight modules when $\ag$ is the untwisted affine Lie algebra of a
simple Lie algebra $\g$ of type $ADE$. The notation will be explained
in the next section.

Let $\lambda$ be a dominant integral weight of $\g$. We write $\lambda
= \sum_{i\in I} m_i\varpi_i$, where $\varpi_i$ is the $i$-th
fundamental weight of $\g$. We consider $\lambda$, $\varpi_i$ as level
$0$ weights of $\ag$ by identifying them with $\lambda - \sum_i m_i
a_i^\vee \Lambda_0$, $\Lambda_i - a_i^\vee \Lambda_0$, where $c =
\sum_i a_i^\vee h_i$ is the central element, and $\Lambda_i$ is the
$i$th fundamental weight of $\ag$.
Let $V(\lambda)$ be the extremal weight module of extremal weight
$\lambda$ with a global crystal base
$(\La(\lambda),\B(\lambda),V^\Z(\lambda))$ 
(see \subsecref{subsec:extremal} for definition).
Let us define a $\Ua$-module
\begin{equation*}
   \widetilde V(\lambda)
   \defeq \bigotimes_{i\in I} V(\varpi_i)^{\otimes m_i},
\end{equation*}
where we take and fix any ordering of $I$ to define the tensor product.
It has $\Uq'(\ag)$-module automorphisms $z_{i,\nu}$ $(i\in I$,
$\nu=1,\dots, m_i$) (see \subsecref{subsec:tensor}).

Set
\(
   \widetilde \La(\lambda)
   \defeq \bigotimes_{i\in I} \La(\varpi_i)^{\otimes m_i},
\)
\(
   \widetilde u_\lambda
   \defeq\bigotimes_{i\in I} u_{\varpi_i}^{\otimes m_i}.
\)
Let $\widetilde \B_0(\lambda)$ be the connected component of the
crystal $\bigotimes_{i\in I} \B(\varpi_i)^{\otimes m_i}$ containing
$\widetilde u_\lambda \bmod q\widetilde\La(\lambda)$.
There is a (subset of) global base $\{ G(b) \mid b\in\B_0(\lambda) \}$
(see \S3.2). Let
\(
   \widetilde \B(\lambda)
   \defeq \{ s(z) b \mid b\in \widetilde \B_0(\lambda) \}
\)
where
\(
  s(z) = \prod_{i\in I} s_{\lambda^{(i)}}(z_{i,1},\dots,
\linebreak[1] z_{i,m_i})
\)
is a product of Schur functions.

There exists a unique $\Ua$-linear homomorphism
\begin{equation*}
   \Phi_\lambda\colon V(\lambda) \to \widetilde V(\lambda)
\end{equation*}
sending $u_\lambda$ to 
\(
   \widetilde u_\lambda
\)
(see \subsecref{subsec:tensor}).

\begin{MainTheorem}\label{thm:main1}
\textup{(1)} $\Phi_\lambda$ is injective.

\textup{(2)} $\Phi_\lambda(\La(\lambda)) \subset \widetilde\La(\lambda)$.

\noindent
Let $\Phi_\lambda^0$ be the induced map
\(
   \La(\lambda)/q\La(\lambda)
   \to \widetilde\La(\lambda)/q\widetilde\La(\lambda).
\)

\textup{(3)} $\Phi_\lambda^0$ gives a bijection between $\B(\lambda)$ and
$\widetilde\B(\lambda)$.

\textup{(4)} $\Phi_\lambda$ maps the global crystal base $\{ G(b) \mid
b\in \B(\lambda) \}$ to $\{ s(z) G(b) \mid b\in\widetilde\B_0(\lambda)\}$.
\end{MainTheorem}

While the author was preparing this article, he learned that Kashiwara
also noticed that his conjecture follows from \cite{BCP} when $\g$ is
of type $ADE$. In fact, some arguments (the proof of the injectivity
of $\Phi_\lambda$, the definition of $(\ ,\ )$, etc.) has been
improved from the original form after the discussion with him.
After the author posted the first version of this paper to the
network archive, he was informed that Jonathan Beck also proved
a part of Kashiwara's conjecture \cite{Beck3}.

\section{Preliminaries}

\subsection{Affine Lie algebra}

Let us fix notations for the untwisted affine Lie algebra $\ag$.
(For a moment we do not assume that $\g$ is of type $ADE$.)
\begin{enumerate}
\item $\widehat I$ : the index set of simple roots,

\item
\(
   \{ \alpha_i\}_{i\in \widehat I}
\)
: the set of simple roots; 
\(
   \{ h_i \}_{i\in \widehat I}
\)
: the set of simple coroots,

\item
\(
   \widehat P^*
   \defeq \bigoplus_{i\in \widehat I} \Z h_i \oplus \Z d
\)
: the dual weight lattice;
$\widehat P = \Hom_\Z(\widehat P^*, \Z)$ : the weight lattice,

\item $\widehat{\mathfrak h} \defeq \widehat P^*\otimes_\Z \Q$ : the
Cartan subalgebra,

\item the simple root $\alpha_i\in\widehat P$ defined by
\(
   \langle h_i, \alpha_j\rangle = a_{ij},
   \langle d, \alpha_j\rangle = \delta_{0j},
\)
where $a_{ij}$ is the Cartan matrix of $\ag$,

\item the fundamental weight $\Lambda_i\in \widehat P$ defined by
\(
   \langle h_i, \Lambda_j \rangle = \delta_{ij},
   \langle d, \Lambda_j\rangle = 0.
\)

\item
\(
   \widehat Q \defeq \bigoplus_{i\in \widehat I} \Z \alpha_i
\) : the root lattice;
\(
   \widehat Q^\vee \defeq \bigoplus_{i\in \widehat I} \Z h_i
\)
: the coroot lattice,

\item
\(
   \widehat Q_+ \defeq \sum_{i\in\widehat I}\Z_{\ge 0}\alpha_i
\);
\(
   \widehat P_+ \defeq \{ \lambda\in\widehat P\mid
   \text{$\langle h_i, \lambda\rangle \ge 0$ for all $i\in\widehat I$}\}
\)
: the set of integral dominant weights,

\item the unique element 
\(
   c = \sum_{i\in\widehat I} a_i^\vee h_i
\)
($a_i^\vee\in\Z_{\ge 0}$) satisfying
\newline
\(
   \left.\left\{ h\in \widehat Q^\vee \,\right|
   \text{$\langle h, \alpha_j\rangle = 0$ for all $j\in\widehat I$}\right\} 
   = \Z c,
\)

\item the unique element
\(
  \delta = \sum_{i\in\widehat I} a_i \alpha_i
\)
($a_i\in\Z_{\ge 0}$) satisfying
\newline
\(
   \left.\left\{ \lambda \in \widehat Q \,\right|
   \text{$\langle h_i, \lambda\rangle = 0$ for all $i\in\widehat I$}\right\} 
   = \Z \delta,
\)

\item the symmetric bilinear form $(\ ,\ )$ on $\widehat{\mathfrak h}^*$,
uniquely characterized by
\(
   \langle h_i, \lambda\rangle =
   \frac{2(\alpha_i, \lambda)}{(\alpha_i, \alpha_i)},
\)
\(
   \langle c, \lambda\rangle = (\delta,\lambda)
\),
for $\lambda\in\widehat{\mathfrak h}^*$,

\item $h \defeq \sum_{i\in\widehat I} a_i$ : the Coxeter number;
$h^\vee \defeq \sum_{i\in\widehat I} a_i^\vee$ : the dual Coxeter
number.
\end{enumerate}

The symmetric bilinear form $(\ ,\ )$ is known to be nondegenerate,
and defines an isomorphism
\(
   \nu\colon\widehat{\mathfrak h}\to\widehat{\mathfrak h}^*
\)
by
\(
   \langle h, \lambda\rangle = (\nu(h), \lambda)
\)
for $\lambda\in\widehat{\mathfrak h}^*$. For example, $\nu(c) =
\delta$. This coincides with one in \cite[\S6]{Kac}.

For $\beta\in\widehat{\mathfrak h}^*$ with $(\beta,\beta)\neq 0$, we
set
\(
   \beta^\vee \defeq \frac{2\beta}{(\beta,\beta)}.
\)
We have
\(
   \nu(h_i) = \alpha_i^\vee.
\)

We have an element $0\in \widehat I$ such that $\{ \alpha_i\mid i\neq
0\}$ is the set of simple roots of $\g$. It is known $a_0^\vee = a_0 = 
1$ for the untwisted affine Lie algebra $\ag$.
We denote $\widehat I\setminus \{0\}$ by $I$.

Let $\operatorname{cl}\colon \widehat{\mathfrak h}^*\to
\widehat{\mathfrak h}^*/\Q\delta$ be the natural projection.
Let
\(
    \widehat{\mathfrak h}^{*0} \defeq
    \{ \lambda\in \widehat{\mathfrak h}^{*0} \mid
    \langle c,\lambda\rangle = 0\},
\)
\(
    \widehat P^0 \defeq \widehat P\cap\widehat{\mathfrak h}^{*0}
\)
(level $0$ weights).
We identify $\operatorname{cl}(\widehat{\mathfrak h}^{*0})\subset
\widehat{\mathfrak h}^*/\Q\delta$ with the dual of the Cartan
subalgebra $\mathfrak h$ of the finite dimensional Lie algebra $\g$,
which is $\bigoplus_{i\in I} \Q h_i$. Similarly we identify
$\operatorname{cl}(\widehat P^0)$ with the weight lattice $P$ of
$\g$. We define the root lattice of $\g$ by
\(
   Q \defeq \bigoplus_{i\in I} \Z\alpha_i
\).
For $i\in I$, we set
\(
   \varpi_i \defeq \Lambda_i - a_i^\vee \Lambda_0 \in \widehat P^0.
\)
Then $\operatorname{cl}(\varpi_i)$ is identified with the $i$th
fundamental weight of $\g$.
Let
\(
  \widehat P^{0,+} \defeq
    \left.\left\{ \lambda \in \widehat P^0\, \right|
    \text{$\langle h_i, \lambda\rangle \ge 0$ for $i\in I$}
  \right\}.
\)
Its projection $\operatorname{cl}(\widehat P^{0,+})$ is the set of
dominant integrable weights of $\g$.
Let
\(
   P^\vee \defeq \Hom_\Z(Q, \Z).
\)
The fundamental coweights $\varpi_i^\vee$ are defined by
\(
   \langle \varpi_i^\vee, \alpha_j\rangle = \delta_{ij}
\)
for $i,j\in I$.
We extend $\varpi_i^\vee$ to a homomorphism $\widehat Q\to\Z$ by
setting $\langle \varpi_i^\vee, \delta\rangle = 0$.

Let $\Delta$ (resp.\ $\Delta_+$) be the set of roots (resp.\ positive
roots) of $\g$. The set of roots $\widehat\aR$ of $\ag$ is given by
$\widehat\aR = \widehat\aR_+\sqcup \widehat\aR_-$, where
\begin{equation*}
  \widehat\aR_+ =
\begin{gathered}
  \{ k\delta + \alpha \mid k\ge 0, \alpha\in\Delta_+\}
  \sqcup
  \{ k\delta \mid k > 0\}
\\
  \sqcup
  \{ k\delta - \alpha \mid k > 0, \alpha\in\Delta_+\},
\end{gathered}
, \qquad
\widehat\aR_- = -\widehat\aR_+.
\end{equation*}
The roots of the form $k\delta \pm \alpha$ ($k\in\Z$,
$\alpha\in\Delta$) are called {\it real\/} roots, while roots
$k\delta$ are called {\it imaginary\/} roots.
The multiplicities of real roots are $1$, and those of imaginary roots
are equal to the rank of $\g$, i.e., $\# I$.

Set
\begin{equation*}
  \begin{gathered}
  \aR_> \defeq \{ k\delta + \alpha \mid  k\ge 0, \alpha\in\Delta^+\},
\quad
  \aR_< \defeq \{ k\delta - \alpha \mid  k > 0, \alpha\in\Delta^+\},
\\
  \aR_0 \defeq \{ k\delta \mid k > 0 \}\times I,
\quad
  \aR \defeq \aR_> \sqcup \aR_0 \sqcup \aR_<.
  \end{gathered}
\end{equation*}
These are sets of roots, counted with multiplicities.

For $i\in \widehat I$, we define the reflection $s_i$ acting on
$\widehat{\mathfrak h}^*$ by
\(
   s_i(\lambda) = \lambda - \langle h_i, \lambda\rangle \alpha_i
\).
Moreover, $s_i$ acts also on $\widehat{\mathfrak h}$ by
\(
   s_i(h) = h - \langle h, \alpha_i\rangle h_i
\).
The actions of $s_i$ preserve $\widehat P$, $\widehat Q$, $\widehat
Q^\vee$ and $\widehat{\mathfrak h}^{*0}$. We have $s_i\delta =
\delta$, $s_i c = c$.
If $i\in I$, the corresponding reflection $s_i$ preserves $\mathfrak
h$, $P$, $P^\vee$ and $Q$. The {\it Weyl group\/} $W$ (resp.\ {\it affine Weyl
group\/} $\widehat W$) of $\g$ (resp.\ $\ag$) is the subgroups of
$\GL(\widehat{\mathfrak h}^*)$ (resp.\ $\GL(\mathfrak h^*)$) generated
by $s_i$ for $i\in \widehat I$ (resp.\ $i\in I$).
We define the {\it extended Weyl group\/} $\eW$ as the semidirect
product
\(
   \eW \defeq W \ltimes P^\vee
\),
using the $W$-action on $P^\vee$.
It is known that $\aW$ is a normal subgroup of $\eW$, and the quotient
$\mathscr T \defeq \eW/\aW$ is a finite group isomorphic to a subgroup
of the group of the diagram automorphisms of $\ag$, i.e., bijections
$\tau\colon I\to I$. Moreover, $\eW$ is isomorphic to $\mathscr
T\ltimes \aW$.

When we consider $\xi\in P^\vee$ as an element of $\eW$, we denote it
by $t_\xi$. We have
\(
   t_\xi(\lambda) = \lambda - \langle \xi, \lambda\rangle \delta
\)
for $\xi\in P^\vee$, $\lambda\in\widehat{\mathfrak h}^{*0}$.

\begin{Lemma}\label{lem:Coxeter}
We have
\begin{equation*}
   \sum_{\alpha\in\widehat\aR_+\cap t_{\varpi_i^\vee}^{-1}(\widehat\aR_-)}
   (\alpha,\xi)
   = h^\vee \langle \varpi_i^\vee, \xi\rangle,
\qquad
   \sum_{\alpha\in\widehat\aR_+\cap t_{\varpi_i^\vee}^{-1}(\widehat\aR_-)}
   (\alpha^\vee,\xi)
   = h \langle \varpi_i^\vee, \xi\rangle.
\end{equation*}
\end{Lemma}

\begin{proof}
From the above description of the root system $\widehat\aR_+$, we have
\begin{equation*}
   \widehat\aR_+\cap t_{\varpi_i^\vee}^{-1}(\widehat\aR_-)
   = \left\{ \beta + n\delta \mid
   \beta\in\Delta_+, 0\le n < \langle \varpi_i^\vee, \beta\rangle \right\}.
\end{equation*}
Therefore
\begin{equation*}
   \sum_{\alpha\in\widehat\aR_+\cap t_{\varpi_i^\vee}^{-1}(\widehat\aR_-)}
   (\alpha,\xi)
   = \sum_{\beta\in\Delta_+}
     (\beta,\xi)\langle \varpi_i^\vee,\beta\rangle
   = \sum_{\beta\in\Delta_+}
     \frac{a_i}{a_i^\vee}(\beta,\xi)(\beta,\varpi_i).
\end{equation*}
We consider the bilinear form on $\mathfrak h^*$ defined by
\begin{equation*}
   \Phi(\xi,\eta) \defeq
    \sum_{\beta\in\Delta_+}(\beta,\xi)(\beta,\eta).
\end{equation*}
From the definition, it is invariant under the Weyl group $W$. So
there is a constant $c$ such that
\(
   \Phi(\xi,\eta) = c(\xi,\eta).
\)
Let $\theta = \delta-\alpha_0$ be the highest root of $\g$. Then we have
\begin{equation*}
   (\theta, \theta)
   = (\alpha_0,\alpha_0) = 2.
\end{equation*}
On the other hand, we have
\begin{equation*}
   \Phi(\theta, \theta)
   = \sum_{\beta\in\Delta_+}(\beta,\theta)(\beta,\theta).
\end{equation*}
If $\beta = \sum_i n_i\alpha_i\in \Delta_+$, we have $0\le n_i\le
a_i$. So we have
\begin{equation*}
\begin{split}
   & (\beta,\theta) = - \sum_i n_i (\alpha_i, \alpha_0) > 0,
\\
   & (\beta,\theta) 
   = (\theta, \theta)
   - \sum_i (n_i-a_i) (\alpha_i,\alpha_0) \le 2,
\end{split}
\end{equation*}
where the equality holds when $\beta = \theta$. (Note that
\(
   (\alpha_i,\alpha_0) = a_{0i}
\)
is a negative integer.) Therefore
\begin{equation*}
\begin{split}
   & \Phi(\theta, \theta) = \sum_{\beta\in\Delta_+} (\beta,\theta) + 2
   = 2(\rho, \theta) + 2
\\
   =\; & 2 \sum_{i\in I} (\varpi_i, \theta) + 2
   = 2 \sum_{i\in I} a_i^\vee + 2 = 2h^\vee,
\end{split}
\end{equation*}
where $\rho$ is the half sum of the positive roots of $\mathfrak h$,
which is known to be equal to $\sum_{i\in I} \varpi_i$. Therefore we
have $c = h^\vee$ and get the first equation. A similar calculation
shows the second equation.
\end{proof}

\subsection{Quantum affine algebra}


Let $\Ua$ be the quantum affine algebra. We follow the notation in
\cite{AK,Kas00}.
We choose a positive integer $d$ such that $(\alpha_i,\alpha_i)/2\in
\Z d^{-1}$ for any $i\in\widehat I$. We set $q_s = q^{1/d}$.
(Later we assume $\g$ is of type $ADE$. Then $d = 1$ and $q_s = q$.)
Then the quantum affine algebra $\Ua$ is the associative algebra over
$\Q(q_s)$ with $1$ generated by elements $e_i$, $f_i$ ($i\in\widehat
I$), $q^h$ ($h\in d^{-1}\widehat P^*$), $q^{\pm c/2}$ with certain
defining relations.
As customary, we set $q_i = q^{(\alpha_i,\alpha_i)/2}$, $t_i =
q^{(\alpha_i,\alpha_i)h_i/2}$, $e_i^{(p)} = e_i^p/[p]_{q_i}!$,
$f_i^{(p)} = f_i^p/[p]_{q_i}!$.

Let $\Uq'(\ag)$ be the quantized enveloping algebra with
$\operatorname{cl}(\widehat P)$ as a weight lattice. It is the
subalgebra of $\Ua$ generated by $e_i$, $f_i$ ($i\in\widehat I$),
$q^h$ ($h\in d^{-1}\bigoplus_i \Z h_i$), $q^{\pm c/2}$.
The quotient $\Uq'(\ag)/(q^{\pm c/2} - 1)$ is denoted by $\Ul$ and
called a {\it quantum loop algebra\/} in \cite{Nak00, Nak01}.

Let $\Ua^+$ (resp.\ $\Ua^-$) be the $\Q(q_s)$-subalgebra of $\Ua$
generated by elements $e_i$'s (resp.\ $f_i$'s).  Let $\Ua^0$ be the
$\Q(q_s)$-subalgebra generated by elements $q^h$ ($h\in d^{-1}\widehat
P^*$). We have the triangular decomposition $\Ua \cong \Ua^+\otimes
\Ua^0 \otimes \Ua^-$.

For $\xi\in \widehat Q$, we define the {\it root space\/} $\Ua_\xi$ by 
\begin{equation*}
   \Ua_\xi \defeq \{ x\in \Ua \mid
     \text{$q^h x q^{-h} = q^{\langle h, \xi\rangle} x$ for all $h\in
     \widehat P^*$}\}.
\end{equation*}

Let $\Ui$ be the $\Z[q_s,q_s^{-1}]$-subalgebra of $\Ua$ generated by
elements $e_i^{(n)}$, $f_i^{(n)}$, $q^h$ for $i\in I$, $n\in\Z_{> 0}$,
$h\in d^{-1}\widehat P^*$.

Let us introduce a $\Q(q_s)$-algebra involutive automorphism $\vee$
and $\Q(q_s)$-algebra involutive anti-automorphisms $*$ and $\psi$ of
$\Ua$ by
\begin{equation*}
\begin{gathered}
   e_i^\vee = f_i, \quad
   f_i^\vee = e_i, \quad
   (q^h)^\vee = q^{-h},
\\
   e_i^* = e_i, \quad
   f_i^* = f_i, \quad
   (q^h)^* = q^{-h},
\\
   \psi(e_i) = q_i^{-1} t_i^{-1} f_i, \quad
   \psi(f_i) = q_i^{-1} t_i e_i, \quad
   \psi(q^h) = q^h.
\end{gathered}
\end{equation*}
We define a $\Q$-algebra involutive automorphism
$\setbox5=\hbox{A}\overline{\rule{0mm}{\ht5}\hspace*{\wd5}}\,$ of $\Ua$
by
\begin{equation*}
\begin{gathered}
   \overline{e_i} = e_i, \quad \overline{f_i} = f_i, \quad
   \overline{q^h} = q^{-h},
\\
   \overline{a(q_s)u} = a(q_s^{-1})\overline{u} \quad
   \text{for $a(q_s)\in\Q(q_s)$ and $u\in\Ua$}.
\end{gathered}
\end{equation*}

In this article, we take the coproduct $\Delta$ on $\Ua$ given by
\begin{equation}\label{eq:comul}
\begin{gathered}
   \Delta q^h = q^h \otimes q^h, \quad
   \Delta e_i = e_i\otimes t_i^{-1} + 1 \otimes e_i,
\\
   \Delta f_i = f_i\otimes 1 + t_i \otimes f_i.
\end{gathered}
\end{equation}

Let us denote by $\varOmega$ the $\Q$-algebra anti-automorphism 
\(
   * \circ \setbox5=\hbox{A}\overline{\rule{0mm}{\ht5}\hspace*{\wd5}}\,
   \circ \vee
\)
of $\Ua$. We have
\begin{equation*}
   \varOmega(e_i) = f_i, \quad \varOmega(f_i) = e_i, \quad
   \varOmega(q^h) = q^{-h}, \quad \varOmega(q_s) = q_s^{-1}.
\end{equation*}

A $\Ua$-module $M$ is called {\it integrable\/} if
\begin{enumerate}
\item all $e_i$, $f_i$ ($i\in I$) are locally nilpotent, and
\item it admits a {\it weight space decomposition\/}:
\[
   M = \bigoplus_{\lambda\in P} M_\lambda, \quad
   \text{where $M_\lambda = \{ u\in M\mid \text{$q^h u
     = q^{\langle h,\lambda\rangle} u$ for all $h\in \widehat P^*$}\}$}.
\]
\end{enumerate}

Let $\Um$ be the modified enveloping algebra \cite[Part IV]{Lu-Book}.
It is defined as 
\begin{equation*}
    \Um \defeq \bigoplus_{\lambda\in \widehat P} \Ua a_\lambda,
\quad
    \Ua a_\lambda \defeq \Ua \left/
    \sum_{h\in \widehat P^*} \Ua (q^h - q^{\langle h, \lambda\rangle})
    \right..
\end{equation*}
Here the multiplication is given by
\begin{equation*}
   a_\lambda x = x a_{\lambda-\xi} \quad
   \text{for $\xi\in\Ua_\xi$},
\qquad
   a_\lambda a_\mu = \delta_{\lambda\mu} a_\lambda,
\end{equation*}
where $a_\lambda$ is considered as the image of $1$ in the above
definition of $\Ua a_\lambda$.

Let $\lambda$, $\mu\in\widehat P_+$. Let $V(\lambda)$ (resp.\
$V(-\mu)$) be the irreducible highest (resp.\ lowest) weight module of 
weight $\lambda$ (resp.\ $-\mu$) \cite[\S3.5]{Lu-Book}. Then there is
a surjective homomorphism
\begin{equation}\label{eq:modify}
   \Ua a_{\lambda-\mu} \ni u \longmapsto u( u_\lambda\otimes u_{-\mu})
   \in V(\lambda)\otimes V(-\mu),
\end{equation}
where $u_\lambda$ (resp.\ $u_{-\mu}$) is a highest (resp.\ lowest)
weight vector of $V(\lambda)$ (resp.\ $V(-\mu)$).

\subsection{Braid group action}

For each $w\in\aW$, there exists an $\Q(q)$-algebra automorphism $T_w$
\cite[\S39]{Lu-Book} (denoted there by $T_{w,1}''$).
%
Also, for any integrable $\Ua$-module $M$, there exists $\Q(q)$-linear
map $T_w\colon M\to M$ satisfying $T_w(xu) = T_w(x) T_w(u)$ for $x\in
\Ua$, $u\in M$ \cite[\S5]{Lu-Book}.
We denote $T_{s_i}$ by $T_i$ hereafter.
By \cite[39.4.5]{Lu-Book} we have
\begin{equation}\label{eq:Omega}
   \varOmega \circ T_w \circ \varOmega = T_w.
\end{equation}

\begin{Lemma}\label{lem:Weyl}
We have
\[
   \left(\psi\circ T_w \circ\psi\right)(x) = (-1)^{N^\vee} q^{-N}\,
   T_{w^{-1}}^{-1}(x) \quad\text{for all $w\in\aW$, $x\in\Ua_\xi$},
\]
where
\[
  N =
  \sum_{\alpha\in\widehat\aR_+\cap w^{-1}(\widehat\aR_-)} (\alpha,\xi)
,\qquad
  N^\vee =
  \sum_{\alpha\in\widehat\aR_+\cap w^{-1}(\widehat\aR_0)} (\alpha^\vee,\xi)
.
\]
\end{Lemma}

\begin{proof}{}
Let $T_{i,-1}''$ be the automorphism defined in
\cite[\S37]{Lu-Book}. A direct calculation shows $\psi\circ T_i
\circ\psi = T_{i,-1}''$. By [loc.\ cit., 37.2.4] we have 
\(
   T_{i,-1}''(x) 
   = (-1)^{\langle h_i, \xi\rangle} q^{-(\alpha_i,\xi)} T_i^{-1}(x)
\)
for $x\in\Ua_\xi$.
Let $w = s_{i_m}\dots s_{i_1}$ be a reduced expression of $w$. Then
\begin{equation*}
   \left(\psi\circ T_w \circ\psi\right)(x)
   =
  (-1)^{N^\vee} q^{-N}
   \left(T_{i_m}^{-1}\dots T_{i_1}^{-1}\right)(x),
\end{equation*}
where
\begin{equation*}
\begin{split}
   N^\vee & = \langle h_{i_1} + s_{i_1}h_{i_2} + \dots +
   s_{i_1}\dots s_{i_{m-1}} h_{i_m},\xi\rangle,
\\
   N &= (\alpha_{i_1} + s_{i_1}\alpha_{i_2} + \dots +
    s_{i_1}\dots s_{i_{m-1}} \alpha_{i_m},\xi).
\end{split}
\end{equation*}
Since we have
\(
   \widehat\aR_+\cap w^{-1}(\widehat\aR_-)
   = 
   \left.\left\{ s_{i_1}\cdots s_{i_{k-1}} \alpha_{i_k}\,\right|\,
     k=1,\dots, m\right\},
\)
we are done.
\end{proof}

As in \cite{Beck,BCP}, the definition of the automorphism $T_w$ of
$\Ua$ can be extended to the case $w\in\eW$ by setting
\[
   \tau e_i = e_{\tau(i)}, \quad
   \tau f_i = f_{\tau(i)}, \quad
   \tau q^{h_i} = q^{h_{\tau(i)}}, \quad
   \tau q^d = q^d.
\]

\subsection{Crystal base}

We shall briefly recall the notion of crystal bases. For the notion of 
(abstract) crystals, we refer to \cite{Kas94,AK}.

For $n\in \Z$ and $i\in \widehat I$, let us define an operator acting
on any integrable $\Ua$-module $M$ by
\begin{gather*}
   \widetilde F_i^{(n)} \defeq
   \sum_{k\ge \max(0,-n)} f_i^{(n+k)} e_i^{(k)} a^n_k(t_i),
\\
   \text{where}\quad
   a^n_k(t_i) \defeq (-1)^k q_i^{k(1-n)}
    \prod_{\nu=1}^{k-1} (1 - q_i^{n+2\nu}).
\end{gather*}
And we set $\te_i \defeq F_i^{(-1)}$, $\tf_i \defeq F_i^{(1)}$.

These operators are different from those used for the definition of
crystal bases in \cite{Kas91}, but it gives us the same crystal bases
by \cite[Proposition~6.1]{Kas00}.

A direct calculation shows
\begin{equation}\label{eq:psi_te}
   \psi(\te_i) = \tf_i.
\end{equation}

Let
\(
  \A_0 \defeq \{ f(q_s)\in\Q(q_s) \mid \text{$f$ is regular at $q_s=0$}\}.
\)

\begin{Definition}
Let $M$ be an integrable $\Ua$-module. A pair $(\La,\B)$ is called a
{\it crystal base\/} of $M$ if it satisfies
\begin{enumerate}
\item $\La$ is a free $\A_0$-submodule of $M$ such that $M\cong
\Q(q_s)\otimes_{\A_0}\La$,
\item $\La = \bigoplus_{\lambda\in P} \La_\lambda$ where
$\La_\lambda = \La\cap M_\lambda$ for $\lambda\in P$,
\item $\B$ is a $\Q$-basis of $\La/q\La \cong \Q\otimes_{\A_0}\La$,
\item $\te_i\La\subset\La$, $\tf_i\La\subset\La$ for all $i\in\widehat
I$,
\item if we denote operators on $\La/q\La$ induced by $\te_i$, $\tf_i$ 
by the same symbols, we have
$\te_i\B\subset\B\sqcup\{0\}$, $\tf_i\B\subset\B\sqcup\{0\}$,
\item for any $b,b'\in\B$ and $i\in \widehat I$, we have
$b' = \tf_i b$ if and only if $b = \te_i b'$.
\end{enumerate}
\end{Definition}

We define functions $\varepsilon_i, \varphi_i\colon\B\to \Z_{\ge 0}$
by
\(
   \varepsilon_i(b) \defeq \max \{ n\ge 0 \mid \te_i^n b\neq 0\},
\)
\(
   \varphi_i(b) \defeq \max \{n \ge 0 \mid \tf_i^n b\neq 0\}.
\)
We set
\(
   \te_i^{\max} b \defeq \te_i^{\varepsilon_i(b)} b,
\)
\(
   \tf_i^{\max} b \defeq \tf_i^{\varphi_i(b)} b.
\)

Let $\setbox5=\hbox{A}\overline{\rule{0mm}{\ht5}\hspace*{\wd5}}\,$ be an
automorphism of $\Q(q_s)$ sending $q_s$ to $q_s^{-1}$. Let $\overline{\A_0}$
be the image of $\A_0$ under
$\setbox5=\hbox{A}\overline{\rule{0mm}{\ht5}\hspace*{\wd5}}\,$, that is,
the subring of $\Q(q_s)$ consisting of rational functions regular at
$q_s=\infty$.

\begin{Definition}
Let $M$ be an integrable $\Ua$-module with a crystale base $(\La,\B)$.
Let $\setbox5=\hbox{A}\overline{\rule{0mm}{\ht5}\hspace*{\wd5}}\,$ be an 
involution of an integrable $\Ua$-module $M$ satisfying
\(
   \overline{xu} = \overline{x}\,\overline{u}
\)
for any $x\in\Ua$, $u\in M$.
Let $M^\Z$ be a $\Ui$-submodule of $M$ such that $\overline{M^\Z} =
M^\Z$, $u - \overline{u} \in (q_s-1)M^\Z$ for $u\in M^\Z$.
We say that {\it $M$ has a global base\/} $(\La,\B,M^\Z)$ if the
following conditions are satisfied
\begin{enumerate}
\item $M \cong \Q(q_s)\otimes_{\Z[q_s,q_s^{-1}]} M^\Z \cong
\Q(q_s)\otimes_{\A_0}\La \cong
\Q(q_s)\otimes_{\overline{\A_0}}\overline{\La}$,
\item $\La\cap\overline{\La}\cap M^\Z \to \La/q_s\La$ is an isomorphism.
\end{enumerate}
As a consequence of the definition, natural homomorphisms
\begin{equation*}
\begin{gathered}
  \A_0\otimes_{\Z}\left(\La\cap\overline{\La}\cap M^\Z\right) \to \La,
\quad
 \overline{\A_0}\otimes_{\Z}\left(\La\cap\overline{\La}\cap
  M^\Z\right)\to \overline{\La},
\\
  \Z[q_s,q_s^{-1}]\otimes_{\Z}\left(\La\cap\overline{\La}\cap
  M^\Z\right) \to M^\Z,
\end{gathered}
\end{equation*}
are isomorphisms.

Let $G$ be the inverse isomorphism
$\La/q_s\La\to\La\cap\overline{\La}\cap M^\Z$. Then $\{ G(b) \mid
b\in\B\}$ is a base of $M$. It is called a {\it global crystal base\/}
of $M$. The above conditions imply $\overline{G(b)} = G(b)$.
\end{Definition}

For a dominant weight $\lambda\in\widehat P_+$, the irreducible
highest weight module $V(\lambda)$ has a global crystal base
\cite{Kas91}. If $\lambda, \mu\in\widehat P_+$, then the tensor
product $V(\lambda)\otimes V(-\mu)$ also has a global crystal base.
Moreover, $\Um$ has a global crystal base 
\( 
  \left( \La(\Um), \B(\Um), \Tilde{\mathbf U}_q^\Z(\ag))\right)
\)
such that the homomorphism \eqref{eq:modify} maps a global base of
$\Um$ to the union of that of $V(\lambda)\otimes V(-\mu)$ and $0$
\cite[Part IV]{Lu-Book}.
Furthermore, the global base is invariant under $*$
\cite[4.3.2]{Kas94}.

\subsection{Extremal vectors}\label{subsec:extremal}

A crystal $\B$ over $\Ua$ is called {\it regular\/} if, for any
$J\subsetneqq\widehat I$, $\B$ is isomorphic (as a crystal over
$\Uq(\g_J)$) to the crystal associated with an integrable
$\Uq(\g_J)$-module. (It was called {\it normal\/} in \cite{Kas94}.)
Here $\Uq(\g_J)$ is the subalgebra generated by $e_j$, $f_j$ ($j\in
J$), $q^h$ ($h\in d^{-1}P^*$).
By \cite{Kas94}, the affine Weyl group $\aW$ acts on any regular
crystal. The action $S$ is given by
\begin{equation*}
   S_{s_i} b =
   \begin{cases}
     \tf_i^{\langle h_i, \wt b\rangle} b
       & \text{if $\langle h_i, \wt b\rangle \ge 0$}, \\
     \te_i^{-\langle h_i, \wt b\rangle} b
       & \text{if $\langle h_i, \wt b\rangle \le 0$}
   \end{cases}
\end{equation*}
for the simple reflection $s_i$. We denote $S_{s_i}$ by $S_i$ hereafter.

\begin{Definition}
Let $M$ be an integrable $\Ua$-module. A vector $u\in M$ with weight
$\lambda\in P$ is called {\it extremal\/}, if the following holds for
all $w\in \aW$:
\begin{equation}\label{eq:extremal}
  \begin{cases}
     e_i T_w u = 0 & \text{if $\langle h_i, w\lambda\rangle \ge 0$},
\\
     f_i T_w u = 0 & \text{if $\langle h_i, w\lambda\rangle \le 0$}.
  \end{cases}
\end{equation}
In this case, we define $S_w u$ so that
\begin{equation*}
    S_i S_w u = 
    \begin{cases}
       f_i^{\left(\langle h_i, w\lambda\rangle\right)} S_w u
       & \text{if $\langle h_i, w\lambda\rangle \ge 0$},
       \\
       e_i^{\left(-\langle h_i, w\lambda\rangle\right)} S_w u
       & \text{if $\langle h_i, w\lambda\rangle \le 0$}.
    \end{cases}
\end{equation*}
This is well-defined, i.e., $S_w u$ depends only on $w$.

Similarly, for a vector $b$ of a regular crystal $B$ with weight
$\lambda$, we say that $b$ is {\it extremal\/} if it satisfies
\begin{equation*}
  \begin{cases}
     \te_i S_w b = 0 & \text{if $\langle h_i, w\lambda\rangle \ge 0$},
\\
     \tf_i S_w b = 0 & \text{if $\langle h_i, w\lambda\rangle \le 0$}.
  \end{cases}
\end{equation*}
\end{Definition}

\begin{Lemma}\label{lem:WeylExtremal}
Suppose that an integrable $\Ua$-module $M$ has a crystal base
$(\La,\B)$. If $u\in \La\subset M$ is an extremal vector of weight
$\lambda$ satisfying $b\defeq u\bmod q\La \in \B$, then $b$ is an
extremal vector, and we have
\begin{equation*}
  S_w u = (-1)^{N^\vee_+} q^{-N_+} T_w u, \quad
  S_w b = S_w u \bmod q\La \quad\text{for all $w\in \aW$},
\end{equation*}
where
\(
\displaystyle
  N_+ =
  \sum_{\alpha\in\widehat\aR_+\cap w^{-1}(\widehat\aR_-)}
  \max\left(\left(\alpha, \lambda\right), 0\right),
\)
and $N^\vee_+$ is given by replacing $\alpha$ by $\alpha^\vee$.
\end{Lemma}

\begin{proof}{}
The equation $S_w b = S_w u\bmod q\La$ follows from the definition of
$S_w$.

If $v\in M_\xi$ satisfies $e_i v = 0$ (resp.\ $f_i v = 0$), we have
\begin{equation*}
   T_i v
   = (-q_i)^{\xi_i} f_i^{(\xi_i)} v
   \qquad
   \left(\text{resp.\ }
   T_i v
   = e_i^{(\xi_i)}v \right),
\end{equation*}
where $\xi_i = \langle h_i, \xi\rangle$.
The rest of the proof is the same as that of \lemref{lem:Weyl}.
\end{proof}

The following follows from a formula for the crystal $\B(\Uqm)$
(see \cite[App.~B]{Kas00}):
\begin{Lemma}\label{lem:extremal}
Let $\lambda\in P^0$. The followings hold for 
\(
   b = b_1\otimes t_\lambda\otimes u_{-\infty}
   \in \B(\Ua a_\lambda)
       = \B(\infty)\otimes T_\lambda\otimes \B(-\infty)
\)
with $\wt b_1\in \Z\delta$\rom:
  
$\te_i b = 0$ or $\tf_i b = 0$ if and only if
$\ve_i(b_1)\le \max(-\langle h_i,\lambda\rangle,0)$.

\end{Lemma}

For $\lambda\in P$, Kashiwara defined the $\Ua$-module $V(\lambda)$
generated by $u_\lambda$ with the defining relation that $u_\lambda$
is an extremal vector of weight $\lambda$ \cite{Kas94}\footnote{He
denoted it by $V^{\max}(\lambda)$.}. It is written as
\begin{equation*}
   V(\lambda) = \Ua a_\lambda/I_\lambda, \qquad
   I_\lambda 
   \defeq \bigoplus_{b\in \B(\Ua a_\lambda)\setminus \B(\lambda)} 
  \Q(q) G(b),
\end{equation*}
where
\(
   \B(\lambda) \defeq
   \{ b\in \B(\Ua a_\lambda)\mid \text{$b^*$ is extremal}\}.
\)
Thus $V(\lambda)$ has a crystal base $(\La(\lambda),\B(\lambda))$
together with a $\Ui$-submodule $V^\Z(\lambda)$ with a global crystal
base, naturally induced from that of $\Ua a_\lambda$.
If $\lambda$ is dominant or anti-dominant, then $V(\lambda)$ is
isomorphic to the highest weight module or the lowest weight
module. So there is no fear of the confusion of the notation.

\subsection{Drinfeld realization}

The quantum affine algebra $\Ua$ has another realization, due to
\cite{Drinfeld,Beck}. It is isomorphic to an associative algebra over
$\Q(q_s)$ with generators $x_{i,r}^\pm$ ($i\in I$, $r\in\Z$),
$q^h$ ($h\in d^{-1}\widehat P^*$), $h_{i,m}^\pm$ ($i\in I$,
$m\in\Z\setminus\{0\}$) with certain defining relations (see
\cite[\S4]{Beck}). The isomorphism depends on the choice of
$o\colon I\to \{\pm 1\}$, and is given by
\begin{gather*}
   x_{i,r}^+ = o(i)^r T_{\varpi_i^\vee}^{-r}(e_i),
 \quad
   x_{i,r}^- = o(i)^r T_{\varpi_i^\vee}^{r}(f_i),
\\
  \left[ x_{i,r}^+, x_{j,s}^- \right]
  = \delta_{ij}\frac{q^{(r-s)c/2} \psi_{i,r+s}^+ - 
  q^{-(r-s)c/2} \psi_{i,r+s}^-}{q - q^{-1}},
\\
  \text{where \ }
  \psi^{\pm}_i(u) \equiv
  \sum_{r=0}^\infty \psi^{\pm}_{i,\pm r} u^{\pm r}
  \defeq t_i^{\pm 1}
   \exp\left(\pm (q_i-q_i^{-1})\sum_{m=1}^\infty h_{i,\pm m} 
     u^{\pm m}\right).
\end{gather*}

By \eqref{eq:Omega} we have
\begin{equation*}
   \varOmega(x^\pm_{i,r}) = x^\mp_{i,-r}, \quad
   \varOmega(h_{i,m}) = h_{i,-m}\quad
   \text{for $i\in I$, $r\in\Z$, $m\in\Z\setminus\{0\}$}.
\end{equation*}

\subsection{The crystal base of $\Ua^+$}\label{subsec:BCP}

Let us recall results in \cite{BCP}.
We assume $\g$ is of type $ADE$ hereafter.
We choose a reduced expression $s_{i_1}\cdots s_{i_N}$ of $2\rho =
2\sum_{i\in I} \varpi_i$ in a suitable way (see [loc.\ cit.] for
detail), and consider a periodic doubly infinite sequence $(\dots,
i_{-1}, i_0, i_1, \dots)$ of $\widehat I$ by setting $i_k = i_{k\bmod
N}$.
Let
\begin{equation*}
   \beta_k \defeq
   \begin{cases}
      s_{i_0}s_{i_{-1}}\cdots s_{i_{k+1}}(\alpha_{i_k})
      & \text{if $k\le 0$},
      \\
      s_{i_1}s_{i_2}\cdots s_{i_{k-1}}(\alpha_{i_k})
      & \text{if $k > 0$}.
   \end{cases}
\end{equation*}
We have
\begin{equation}\label{eq:sign}
   \aR_> = \{ \beta_k \mid k \le 0\}, \quad
   \aR_< = \{ \beta_k \mid k > 0 \}.
\end{equation}

We define
\begin{equation*}
   E_{\beta_k}^{(n)} \defeq
   \begin{cases}
     T_{i_0}^{-1} T_{i_{-1}}^{-1}\dots T_{i_{k+1}}^{-1}(e_{i_k}^{(n)})
     & \text{if $k\le 0$},
     \\
     T_{i_1} T_{i_2} \dots T_{i_{k-1}}(e_{i_k}^{(n)})
     & \text{if $k > 0$}.
   \end{cases}
\end{equation*}
We denote $E_{\beta_k}^{(1)}$ by $E_{\beta_k}$. These are {\it root
vectors\/} for $\aR_>$ and $\aR_<$.
By \cite[40.1.3]{Lu-Book} we have $E_{\beta_k}^{(n)}\in \Ua^+$.
Explicit relations among $E_{\beta_k}^{(n)}$ and $x_{i,r}^\pm$ can be
found in \cite[Lemma~1.5]{BCP}.


We define $P_{m,i}$ ($m > 0$, $i\in I$) by
\begin{equation*}
   1 + \sum_{m > 0} P_{m,i} u^m
   = \exp\left(
   - \sum_{m > 0} \frac{ (o(i) q^{c/2} u )^r h_{i,r}}{[r]_{q_i}} \right).
\end{equation*}
We also define $\widetilde P_{m,i}\in \Ua^+$ by replacing $h_{i,r}$ by 
$-h_{i,r}$. These are root vectors for $\aR_0$. We also set
\(
   P_{-m,i} = \varOmega(P_{m,i})
\)
\(
  (m > 0, i\in I).
\)

Let $\bc \colon \aR\to\Z_{\ge 0}$ be a map such that $\bc(\alpha) =
0$ except for finitely many $\alpha$. We denote its restrictions to
$\aR_>$, $\aR_>$, $\aR_0$ by $\bc_>$, $\bc_<$, $\bc_0$ respectively.
We define $E_{\bc_>}, E_{\bc_<} \in \Ua^+$ by
\begin{equation*}
   E_{\bc_>} \defeq E_{\beta_0}^{(\bc(\beta_0))}
   E_{\beta_{-1}}^{(\bc(\beta_{-1}))}\cdots,
\qquad
   E_{\bc_<} \defeq \cdots E_{\beta_2}^{(\bc(\beta_2))}
   E_{\beta_1}^{(\bc(\beta_1))}.
\end{equation*}
Next, given $\bc_0$, we associate an $I$-tuple of partitons
$(\lambda^{(i)})_{i\in I}$ as
\begin{equation*}
    \lambda^{(i)} \defeq (1^{\bc_0(\delta,i)} 2^{\bc_0(2\delta,i)} \cdots
    k^{\bc_0(k\delta,i)}\cdots).
\end{equation*}
As in \cite{Mac} we denote it also in another notation:
\begin{equation*}
    \lambda^{(i)} = (\lambda^{(i)}_1, \lambda^{(i)}_2, \dots).
\end{equation*}
We define the corresponding Schur function
\begin{equation*}
   S_{\bc_0} \defeq
   \prod_{i\in I}
   \det\left(\widetilde P_{\lambda^{(i)}_k-k+l,i}\right)_{1\le k,l\le t},
\end{equation*}
where $t\ge l(\lambda^{(i)})$. Note that $\widetilde P_{m,i}$
corresponds to a complete symmetric function, while $P_{m,i}$
corresponds to an elementary symmetric function, up to sign.

Now a main result of \cite{BCP} says that
\begin{enumerate}
\item $B_{\bc} \defeq \overline{E_{\bc_>} \cdot S_{\bc_0} \cdot
E_{\bc_<}}$ is contained in $\La(\infty)\cap \Ui^+$,
\item $\{ B_{\bc} \bmod q\La(\infty) \mid \bc\in \Z_{\ge 0}^{\aR} \}$
is the crystal base of $\Ua^+$.
\end{enumerate}

Set
\(
   (\Z_{\ge 0}^{\aR_0})(\lambda)
   \defeq
   \left.\left\{ \bc_0\in\Z_{\ge 0}^{\aR_0}\, \right|\, 
   \text{$l(\lambda^{(i)}) \le \langle h_i,\lambda\rangle$ for all $i\in I$}
   \right\},
\)
where $(\lambda^{(i)})_{i\in I}$ is the $I$-tuple of partition
corresponding to $\bc_0$ as above.

We apply $\vee$ to the above crystal base to get
\begin{equation*}
   F_{\bc_>} \defeq (E_{\bc_>})^\vee, \quad
   F_{\bc_<} \defeq (E_{\bc_<})^\vee, \quad
   S_{\bc_0}^- \defeq (S_{\bc_0})^\vee.
\end{equation*}

\subsection{Extremal weight modules and the Drinfeld realization}

Extremal weight modules are defined in terms of Chevalley
generators. We shall rewrite the definition in terms of Drinfeld
generators, and derive several easy consequences in this subsection.

The following is a consequence of \cite[Theorem~5.3]{Kas00}.

\begin{Lemma}\label{lem:l-dom}
Let $u$ be a vector of an integrable $\Uq'(\ag)$-module $M$ with
weight $\lambda\in \widehat P^{0,+}$. Then the following conditions
are equivalent\rom:

\textup{(1)} $u$ is an extremal vector.

\textup{(2)}
\(
   x^+_{i,r} u = 0
\)
for all $i\in I$, $r\in\Z$.
\end{Lemma}




\begin{Remark}
The extremal weight module $V(\lambda)$ is isomorphic to the Weyl
module $W_q(\lambda)$ introduced by Chari-Pressley \cite{CP00}. This
result was refered as `an unpublished work' of Kashiwara in
[loc.\ cit., Proposition~4.5]. Let us give Kashiwara's proof here.
Let $\lambda = \sum_{i\in I} m_i\varpi_i\in \widehat P^{0,+}$.
Then $W_q(\lambda)$ is integrable and contains a vector $w_\lambda$ of
weight $\lambda$ which satisfies the above condition~(2). Therefore,
there is a unique $\Ua$-linear homomorphism $V(\lambda)\to
W_q(\lambda)$, sending $v_\lambda$ to $w_\lambda$.
(The integrablity of $W_q(\lambda)$ was proved via the isomorphism 
$V(\lambda)\cong W_q(\lambda)$ in [loc.\ cit.]. So one must give
another proof of the integrablity as sketched in [loc.\ cit.].)
Since $W_q(\lambda)$ is generated by $w_\lambda$ by definition, the
homomorphism is surjective.
On the other hand, any integrable $\Ua$-module generated by a vector
$u$ of weight $\lambda$ satisfying the above condition~(2) is a
quotient of $W_q(\lambda)$ [loc.\ cit., Proposition~4.6]. Therefore
$V(\lambda)$ and $W_q(\lambda)$ are isomorphic.
\end{Remark}

\begin{Corollary}\label{cor:P}
Let $u$ be an extremal vector with weight $\lambda\in \widehat
P^{0,+}$. Then
\(
   S^-_{\bc_0} u = \overline{S_{\bc_0}^*} u = 0 
\)
   for $\bc_0\notin(\Z_{\ge 0}^{\aR_0})(\lambda)$.
\end{Corollary}

\begin{proof}
It is enough to show the assertion for $u = u_\lambda\in V(\lambda)$.
We have a $\Q(q)$-vector space isomorphism
\begin{equation*}
    V(\lambda) \ni x u_\lambda \mapsto
      x^\vee u_{-\lambda} \in V(-\lambda). 
\end{equation*}
Therefore it is enough to show $S_{\bc_0} u_{-\lambda}
= \Omega(S_{\bc_0}) u_{-\lambda} = 0$.
By \cite[Proposition~4.3]{CP00}, which is applicable thanks to
\lemref{lem:l-dom}, we have
\begin{equation*}
   P_{m,i} u_{-\lambda} = 0
   \qquad\text{for $|m| > \langle h_i, \lambda\rangle$}.
\end{equation*}
(More precisely, we apply [loc.\ cit.] after composing an automorphism 
$x_{i,r}^\pm \mapsto -x_{i,-r}^\mp$, $h_{i,m} \mapsto -h_{i,-m}$.)
Now the assertion follows from a standard result in the theory of
symmetric polynomials.
\end{proof}

\section{A study of extremal weight modules}

\subsection{Fundamental representations}

By \cite[\S5.2]{Kas00} there is a unique $\Uq'(\ag)$-linear
automorphism $z_i$ of $V(\varpi_i)$ with weight $\delta$, which sends 
$u_{\varpi_i}$ to $u_{\varpi_i+\delta}$. (Note that $d_i$ in
\cite[\S5.2]{Kas00} is equal to $1$ for untwisted $\ag$.)

\begin{Proposition}\label{prop:Drinfeld}
\(
   h_{i,1} u_{\varpi_i} = o(i) (-1)^{1-h} q^{-h^\vee} z_i u_{\varpi_i}
\).
\end{Proposition}

\begin{proof}
We have
\begin{equation*}
   h_{i,1} u_{\varpi_i}
   = t_i^{-1} \left[x^+_{i,1}, x^-_{i,0}\right] u_{\varpi_i}
   = o(i) t_i^{-1} T_{\varpi_i^\vee}^{-1} (e_i) f_i u_{\varpi_i}.
\end{equation*}
Let us write $T_{\varpi_i^\vee} = \tau T_w$ with $w\in\aW$. Then
\lemref{lem:WeylExtremal} implies
\begin{equation}\label{eq:temp}
    T_{\varpi_i^\vee}^{-1} (e_i) f_i u_{\varpi_i}
    = 
    (-1)^{N^{\prime\vee}_+} q^{N'_+}
    S_w^{-1} \left(e_{\tau^{-1}(i)} S_w (f_i u_{\varpi_i})\right),
\end{equation}
where
\(
  N'_+ =
  \sum_{\alpha\in\widehat\aR_+\cap w^{-1}(\widehat\aR_-)}
  \max((\alpha,s_i\varpi_i),0)
  - \max((\alpha,\varpi_i), 0)
,
\)
and $N^{\prime\vee}_+$ is given by replacing $\alpha$ by $\alpha^\vee$.
Since
\(
   \widehat\aR_+\cap w^{-1}(\widehat\aR_-) =
   \widehat\aR_+\cap t_{\varpi_i^\vee}^{-1}(\widehat\aR_-) =
   \{ \beta+n\delta \mid \beta\in\Delta_+,
     0\le n < \langle \varpi_i, \alpha\rangle\},
\)
we have
\begin{equation*}
   \max((\alpha, \varpi_i),0) = (\alpha, \varpi_i),
\quad
   \max((\alpha, s_i\varpi_i),0)
   = 
   \begin{cases}
      0 & \text{if $\alpha = \alpha_i$},
      \\
      (\alpha, s_i\varpi_i) & \text{otherwise}.
   \end{cases}
\end{equation*}
Therefore
\begin{equation*}
   N'_+ = (\alpha_i,\varpi_i) - 
   \sum_{\alpha\in\widehat\aR_+\cap w^{-1}(\widehat\aR_-)}
     (\alpha,\alpha_i)
   = (\alpha_i,\varpi_i) - h^\vee,
\end{equation*}
where we have used \lemref{lem:Coxeter}. Similarly we have
\(
   N^{\prime\vee}_+ = 1 - h
\).
Now the assertion follows from the definition of the Weyl group action
$S$.
\end{proof}

\begin{Remark}
Let $W(\varpi_i) \defeq V(\varpi_i)/(z_i - 1)V(\varpi_i)$. This is a
finite dimensional irreducible $\Uq'(\ag)$-module \cite[\S5.2]{Kas00}.
The above proposition says that $W(\varpi_i)$ has the Drinfeld
polynomial
\begin{equation*}
   P_j(u) = 
   \begin{cases}
     1 & \text{if $j\neq i$},
     \\
     1 + o(i)(-1)^h q^{-h^\vee}u & \text{if $j=i$}.
   \end{cases}
\end{equation*}
\end{Remark}

\begin{Proposition}\label{prop:fund}
\(
    (\tilde P_{\pm 1,i})^\vee u_{\varpi_i} = z_i^\pm u_{\varpi_i}
\).
\end{Proposition}

\begin{proof}
Let us endow a new $\Ua$-module structure on $V(-\varpi_i)$ by
\begin{equation*}
   x \cdot u \defeq x^\vee \cdot u,
   \quad(x\in \Ua, u\in V(-\varpi_i)).
\end{equation*}
We denote it by $V(-\varpi_i)^\vee$. Then there is a $\Ua$-module
isomorphism $V(\varpi_i)\cong V(-\varpi_i)^\vee$ sending
$u_{\varpi_i}$ to $u_{-\varpi_i}$. Using this isomorphism, we can
calculate $(\tilde P_{\pm 1,i})^\vee u_{\varpi_i}$ exactly as in the
above proposition (in fact, more easily) to get the assertion.
\end{proof}

\subsection{Tensor product modules}\label{subsec:tensor}

Let $\lambda = \sum_{i\in I} m_i \varpi_i\in \widehat P^{0,+}$. We
define a $\Ua$-module $\widetilde V(\lambda)$, $\widetilde
\La(\lambda)$, $\widetilde \B(\lambda)$, $\widetilde u_\lambda$ as in
the introduction.
Let $z_{i,\nu}$ $(i\in I$, $\nu=1,\dots, m_i$) be the
$\Uq'(\ag)$-linear automorphism of $\widetilde V(\lambda)$ obtained by
the action of $z_i\colon V(\varpi_i)\to V(\varpi_i)$ on the $\nu$-th
factor. Obviously they are commuting: $z_{i,\nu} z_{j,\mu} = z_{j,\mu}
z_{i,\nu}$. 
Let
\begin{equation*}
\begin{gathered}
   \Breve V(\lambda)\defeq
   \Ua[z_{i,\nu}^\pm]_{i\in I,\nu=1,\dots,m_i} \cdot\widetilde u_\lambda,
\quad
   \Breve \La(\lambda)\defeq
   \widetilde \La(\lambda)\cap \Breve V(\lambda),
\\
   \Breve \B(\lambda)\defeq
   \bigotimes_{i\in I} \B(\varpi_i)^{\otimes m_i},
\quad
   \Breve V^\Z(\lambda) \defeq
   \bigotimes_{i\in I} \left(V(\varpi_i)^\Z\right)^{\otimes m_i}\cap
   \Breve V(\lambda).
\end{gathered}
\end{equation*}

By \cite[\S8]{Kas00}, the submodule $\Breve V(\lambda)$ has
\begin{enumerate}
\item the unique bar involution
$\setbox5=\hbox{A}\overline{\rule{0mm}{\ht5}\hspace*{\wd5}}\,$
satisfying
\newline
$\overline{xu} = \overline{x}\, \overline{u}$ for $x\in
\Ua[z_{i,\nu}^\pm]_{i\in I,\nu=1,\dots,m_i}$, $u\in \Breve
V(\lambda)$,
\item the crystal base $(\Breve \La(\lambda),\Breve \B(\lambda))$, and
\item the $\Ui$-submodule $\Breve V^\Z(\lambda)$ and
the global crystal base $\{ G(b) \mid b\in\Breve\B(\lambda)\}$.
\end{enumerate}

The module $\widetilde V(\lambda)$ contains the extremal vector
\(
   \widetilde u_\lambda
\)
of weight $\lambda$. Therefore there exists a unique $\Ua$-linear
homomorphism
\(
   \Phi_\lambda\colon V(\lambda) \to \widetilde V(\lambda)
\)
sending $u_\lambda$ to 
\(
   \widetilde u_\lambda
\).
The image is contained in $\Breve V(\lambda)$.

Recall that a function $\bc_0\in \aR_0\to \Z_{\ge 0}$ defines
an $I$-tuple of partitions $(\lambda^{(i)})_{i\in I}$ as
\subsecref{subsec:BCP}.
We define an endomorphism of $\widetilde V(\lambda)$ by
\[
   s_{\bc_0}(z^\pm) \defeq
   \prod_{i\in I} s_{\lambda^{(i)}}(z_{i,1}^\pm,\dots, z_{i,m_i}^\pm),
\]
where $s_{\lambda^{(i)}}$ is the Schur polynomial corresponding to the
partition $\lambda^{(i)}$. If $l(\lambda^{(i)}) > m_i$, it is
understood us $0$.

\begin{Proposition}\label{prop:tensor}
\(
   \Phi_\lambda(S_{\bc_0}^- u_\lambda)
   = s_{\bc_0}(z)\cdot \widetilde u_\lambda
\)
,
\(
   \Phi_\lambda(\overline{S_{\bc_0}^*} u_\lambda)
   = s_{\bc_0}(z^{-1})\cdot \widetilde u_\lambda
\).
\end{Proposition}

\begin{proof}
On level $0$ modules, we have
\[
   \Delta h_{i,\pm m} = h_{i,\pm m} \otimes 1 + 1\otimes h_{i,\pm m}
   + \text{a nilpotent term}
\]
by \cite{Da}. Up to sign, the transition between $h_{i,m}$'s and
$P_{k,i}$'s is the same as that between power sums and elementary
symmetric functions. The above equation means that $\Delta$ coincides
with the standard coproduct on symmetric polynomials modulo nilpotent
terms \cite[Chap.\ I, \S5, Ex.\ 25]{Mac}. Therefore we have
\[
   \Delta P_{k,i}
   = \sum_{s=0}^{k} P_{s,i}\otimes P_{k-s,i}
   + \text{a nilpotent term}.
\]
Using \corref{cor:P} and \propref{prop:fund}, we have the assertion.
\end{proof}

\subsection{Detemination of extremal vectors}

\begin{Proposition}\label{prop:extremal}
Suppose $\lambda\in\widehat P^{0,+}$.  Consider
\(
   B_{\bc} = \overline{F_{\bc_>}\cdot S_{\bc_0}^- \cdot F_{\bc_<}}
\)
with $\wt B_{\bc}\in \Z\delta$, and set 
\(
   b_1 \defeq B_{\bc} \bmod q\La(\infty) \in \B(\infty)
\)
and
\(
  b \defeq b_1\otimes t_\lambda\otimes u_{-\infty}\in
  \B(\Uqm a_\lambda)
\).
If $b$ and $b^*$ are extremal, then we have
\(
     \bc_> \equiv 0 \equiv \bc_<
\)
and $\bc_0\in(\Z_{\ge 0}^{\aR_0})(\lambda)$.
\end{Proposition}

\begin{proof}
Assume $\bc_>\not\equiv 0$ and take the largest number $k\le 0$
satisfying $\bc(\beta_k) \neq 0$. Let $w = s_{i_{0}} s_{i_{-1}} \cdots 
s_{i_{k+1}}$.

Since $b^*$ is extremal, we can consider $b$ as an element of
$\B(\lambda)$. We have
\begin{equation*}
   b = B_{\bc} u_\lambda \bmod q\La(\lambda).
\end{equation*}
By \lemref{lem:WeylExtremal}, we have
\begin{equation*}
   S_w^{-1} b = (-1)^{N^\vee} q^{N} 
      T_w^{-1}(B_{\bc}) \cdot S_w^{-1}( u_\lambda) \bmod q\La(\lambda)
\end{equation*}
for some integers $N^\vee$, $N$.
By \cite[8.2.2]{Kas94} there exists a $\Ua$-linear isomorphism
\begin{equation*}
    V(\lambda) \to V(w^{-1}\lambda);
    \qquad S_w^{-1}(u_\lambda) \mapsto u_{w^{-1}\lambda},
\end{equation*}
respecting the crystal bases.
Therefore we have
\begin{equation*}
   (-1)^{N^\vee} q^{N} T_w^{-1} (B_{\bc})
   u_{w^{-1}\lambda} \bmod q\La({w^{-1}\lambda})
   \in \B({w^{-1}\lambda}).
\end{equation*}
(In fact, this is equal to $S_{w^{-1}}^* S_{w^{-1}} b$.)
Let us denote this by $b_1'\otimes t_{w^{-1}\lambda}\otimes b_2'$.

We have
\begin{equation*}
   T_w^{-1} (B_{\bc}) 
   = T_w^{-1} (\overline{F_{\bc_>}})\cdot
     T_w^{-1} (\overline{S_{\bc_0}^-})\cdot
     T_w^{-1} (\overline{F_{\bc_<}}).
\end{equation*}
It is clear that
\(
    T_w^{-1} (\overline{F_{\bc_<}})\in \Ua^-\cap T_{i_k}\Ua^-.
\)
We also have
\(
    T_w^{-1} (\overline{S_{\bc_0}^-})\in \Ua^-\cap T_{i_k}\Ua^-
\)
by \cite[Lemma~2]{Beck2}. (More precisely, we apply [loc.\ cit.] after
composing
$\setbox5=\hbox{A}\overline{\rule{0mm}{\ht5}\hspace*{\wd5}}\,\circ
\vee$. Note that \( T_{w}^{-1} =
\setbox5=\hbox{A}\overline{\rule{0mm}{\ht5}\hspace*{\wd5}}\,\circ
\vee\circ T_w\circ
\setbox5=\hbox{A}\overline{\rule{0mm}{\ht5}\hspace*{\wd5}}\,\circ \vee
\) by \cite[39.4.5]{Lu-Book}.)
Moreover, by our choice of $k$, we have
\begin{equation*}
  T_w(\overline{F_{\bc_>}}) = f_{i_k}^{(\bc(\beta_k))}\, 
      T_{i_k} (f_{i_{k-1}}^{(\bc(\beta_{k-1}))}) \cdots
  \in
  f_{i_k}^{(\bc(\beta_k))} \left( \Ua^- \cap T_{i_k}\Ua^-\right).
\end{equation*}
Therefore we have
\begin{equation*}
   b'_2 = u_{-\infty},
\qquad
   b'_1 = T_w^{-1}(B_{\bc}) \bmod q\La(\infty), \qquad
   \ve_{i_k}(b'_1) = \bc(\beta_k),
\end{equation*}
where the last equality follows from \cite[38.1.6]{Lu-Book}.
Since $b_1'\otimes t_{w^{-1}\lambda}\otimes u_{-\infty}$ is extremal,
\lemref{lem:extremal} implies
\begin{equation}\label{eq:tt}
   \bc(\beta_k) \le \max(-\langle h_{i_k}, w^{-1}\lambda\rangle,0).
\end{equation}
However, we have
\(
   \langle h_{i_k}, w^{-1}\lambda \rangle 
   = (w \alpha_{i_k}^\vee, \lambda) \ge 0
\)
for $\lambda\in\widehat P^{0,+}$, because
$w\alpha_{i_k}\in\widehat\aR_>$ by \eqref{eq:sign}.
So the right hand side of \eqref{eq:tt} is $0$, and this contradicts
with the choice of $k$. Therefore $\bc_> \equiv 0$. Applying $*$, we
similarly get $\bc_< \equiv 0$.
Now the last assertion is a consequence of \corref{cor:P}.
\end{proof}

\begin{proof}[Proof of \thmref{thm:main1}]
We first prove (2), (3), (4) and then (1).

(2) Recall that any vector $b\in\B(\lambda)$ is connected to an
extremal vector \cite[9.3.3]{Kas94}. Moreover, an extremal vector can
be mapped by $\tf_i^{\max}$ to an extremal vector of the form
$b_1\otimes t_\lambda\otimes u_{-\infty}$. (See \cite[Proof of
Theorem~5.1]{Kas00}). Therefore
\begin{equation*}
   \B(\lambda) = \left\{ X_l \cdots X_1 \overline{S_{\bc_0}^-}
   \bmod q\La(\Lambda) \left|\,
   \bc_0\in(\Z_{\ge 0}^{\aR_0})(\lambda),
   \ \text{$X_\mu$ is $\te_i$ or $\tf_i$}
   \right\}\right.\setminus \{0\}
\end{equation*}
by \propref{prop:extremal}. Then $\La(\lambda)$ is spanned by $\{ X_l
\cdots X_1 \overline{S_{\bc_0}^-}\}$ over $\Z[q]$, by Nakayama's lemma.
Note that $\Phi_\lambda$ commutes with the operators $\te_i$, $\tf_i$
and $\widetilde \La(\lambda)$ is invariant under $\te_i$, $\tf_i$.
Therefore it is enough to show that $\Phi_\lambda(\overline{S_{\bc_0}^-})\in
\widetilde\La(\lambda)$. But this follows from \propref{prop:tensor}.

(3) By \propref{prop:tensor}, we have
\[
   \Phi_\lambda^0(\overline{S_{\bc_0}^-}\bmod q\La(\lambda))
   \in\widetilde\B(\lambda)
   \qquad\text{for $\bc_0\in(\Z_{\ge 0}^{\aR_0})(\lambda)$}.
\]
As in the proof of (1), we conclude that
$\Phi_\lambda^0(\B(\lambda))\subset\widetilde\B(\lambda)\sqcup\{ 0\}$.
From the definition, it is obvious that the image contains
$\widetilde\B(\lambda)$.
Consider $\Ker\Phi^0_\lambda\cap\B(\lambda)$. It is invariant under
$\te_i$, $\tf_i$. Since any vector is connected to an extremal vector, 
$\Ker\Phi^0_\lambda\cap\B(\lambda)$ contains an extremal vector if it
is nonempty. But we already checked that every extremal vector is
mapped to a nonzero vector. Hence $\Phi^0_\lambda|_{\B(\lambda)}$ is
injective.

(4) By the uniqueness, $\Phi_\lambda$ respects the bar involutions on
$V(\lambda)$ and $\widetilde V(\lambda)$. Since
$V^\Z(\lambda) = \Ui u_\lambda$, we have
\(
   \Phi_\lambda(V^\Z(\lambda))\subset \Breve V^\Z(\lambda)
\).
Therefore we have
\[
  \Phi_\lambda\left(\La(\lambda)\cap\overline{\La(\lambda)}\cap
    V^\Z(\lambda)\right)
  \subset \Breve\La(\lambda)\cap\overline{\Breve\La(\lambda)}\cap
  \Breve V^\Z(\lambda).
\]
Now the assertion follows from (3).

(1) It is easy to see that $\widetilde\B(\lambda)$ is linearly
independent. Therefore $\Phi_\lambda^0\colon
\La(\lambda)/q\La(\lambda)\to
\widetilde\La(\lambda)/q\widetilde\La(\lambda)$ is injective.

Let $\{G(b)\}$ be the global crystal base of $V(\lambda)$. Let $0\neq
\sum f_b(q) G(b)\in\Ker\Phi_\lambda$. Multiplying a power of $q$, we
may assume $f_b(q)\in\A_0$ for all $b$ and $f_{b_0}(0)\neq 0$ for some
$b_0$.  Then $\sum f_b(0) b\in \La(\lambda)/q\La(\lambda)$ is mapped
to $0$ by $\Phi_\lambda^0$. The injectivity of $\Phi_\lambda^0$
implies that $f_b(0) = 0$ for all $b$. This is a contradiction.
\end{proof}

\begin{Remark}\label{rem:exgl}
\thmref{thm:main1} together with \propref{prop:tensor} implies that
$S_{\bc_0}^- u_\lambda$ and $\overline{S_{\bc_0}^*} u_\lambda$ are elements
of the global base.
\end{Remark}

\subsection{Standard modules}\label{subsec:standard}

Let us briefly recall the properties of the universal standard module
$M(\lambda)$ with a weight $\lambda = \sum m_i \varpi_i\in \widehat
P^{0,+}$ introduced in \cite{Nak00,Nak01}. (We do not review its
definition, which is based on quiver varieties.)
Let $G_\lambda \defeq \prod_i \GL_{m_i}(\C)$. Its maximal torus
consisting of diagonal matrices is denoted by $H_\lambda$. Their
representation rings are denoted by $R(G_\lambda)$,
$R(H_\lambda)$ respectively. They are isomrphic to
\(
   \bigotimes_i \Z[x_{i,1}^\pm,\dots, x_{i,m_i}^\pm]^{\mathfrak S_{m_i}}
\)
and
\(
   \bigotimes_i \Z[x_{i,1}^\pm,\dots, x_{i,m_i}^\pm]
\)
respectively.
The universal standard module $M(\lambda)$ is a
$\Uq^{\prime\Z}(\ag)\otimes_{\Z} R(G_\lambda)$-module which is
integrable (in fact, it satisfies a stronger condition `{\it
l\/}--integrability') and contains a vector $[0]_\lambda$ with
\begin{gather*}
   x_{i,r}^+[0]_\lambda = 0\quad\text{for any $i\in I$, $r\in \Z$},
\quad
   q^h [0]_\lambda = q^{\langle h,\lambda\rangle} [0]_\lambda,
\\
   M(\lambda) = 
   \left(\Uq^{\prime\Z}(\ag)\otimes_{\Z} R(G_\lambda)\right)[0]_{\lambda},
\\
   \psi_i^\pm(u)[0]_\lambda
   = q^{m_i} \left( \prod_{\nu=1}^{m_i}
     \frac{1-q^{-1}x_{i,\nu}u}{1-q x_{i,\nu}u}\right)^\pm
   [0]_\lambda,
\end{gather*}
where $(\ )^\pm$ denotes the expansion at $u=0$ and $\infty$
respectively.
(In fact, we have $M(\lambda) = \Uq^{\prime\Z}(\ag)[0]_{\lambda}$ by
the proof of \thmref{thm:main1}.)
Moreover, $M(\lambda)$ is free of finite rank as an
$R(G_\lambda)$-module. And $M(\lambda)$ is simple if we tensor
the quotient field of $\Z[q,q^{-1}]\otimes R(G_\lambda)$.

On the other hand, we have a 
\(
   \bigotimes_{i\in I}
   \Z[z_{i,1}^\pm,\dots, z_{i,m_i}^\pm]^{\mathfrak S_{m_i}}
\)
-module structure on $V(\lambda)$ given by
\(
   s_{\bc_0}(z) u_\lambda = S_{\bc_0}^- u_\lambda
\)
and 
\(
   s_{\bc_0}(z^{-1}) u_\lambda = \overline{S_{\bc_0}^*} u_\lambda
\)
by the above discussion. We make it a
\(
   R(G_\lambda) = \bigotimes_{i\in I}
   \Z[x_{i,1}^\pm,\dots, x_{i,m_i}^\pm]^{\mathfrak S_{m_i}}
\)
-module structure by setting
\(
   x_{i,\nu} = o(i)(-1)^{1-h}q^{-h^\vee}z_{i,\nu}.
\)

\begin{MainTheorem}
There exists a unique
\(
  \Uq^{\prime\Z}(\ag)\otimes_\Z R(G_\lambda)
\)
-isomorphism $V^\Z(\lambda)\to M(\lambda)$ sending $u_\lambda$ to
$[0]_\lambda$.
\end{MainTheorem}

This result follows from \thmref{thm:main1} as explained in
\cite[1.23]{Nak01}. The calculation of Drinfeld polynomial, which was
not given there, is done in \propref{prop:Drinfeld}.

\noindent {\bf Correction to \cite{Nak01}}:

Delete $\mathfrak S_{\lambda_1}\times\cdots\times\mathfrak S_{\lambda_n}$ in
Theorem~1.22.

Replace $R(G_\lambda)$ in page 411, line 5 by $R(H_\lambda)$.

Delete `and forgetting the symmetric group invariance'
in Remark~1.23.

Replace `the submodule above' in line 8, `the submodule
\newline
\(
  \Uli[x_{k,\nu}]_{k\in I,\nu=1,\dots,\lambda_k}
  \bigotimes_{k\in I} [0]_{\Lambda_k}^{\otimes \lambda_k}
\).

\section{A bilinear form}

Kashiwara proved that the crystal base $\B(\lambda)$ is an orthonomal
base with respect to a natural bilinear form when $\lambda$ is
dominant \cite[5.1.1]{Kas91}. We prove a similar result for
$\lambda\in\widehat P^{0,+}$ in this section. This generalizes a
result of Varagnolo-Vasserot \cite[Theorem A]{VV2} from fundamental
representations to arbitray $\lambda$.

\begin{Proposition}[Kashiwara]\label{prop:bilinear}
The extremal weight module $V(\lambda)$ has a unique bilinear form $(\ 
,\ )$ satisfying
\begin{gather}
   (u_\lambda, G(b)) = 
   \begin{cases}
      1 & \text{if $G(b) = u_\lambda$}, \\
      0 & \text{otherwise}
   \end{cases}\label{eq:bil1}
\\
   (x u, v) = (u, \psi(x)v) \quad\text{for $x\in\Ua$, $u,v\in V(\lambda)$}.
   \label{eq:bil2}
\end{gather}
\end{Proposition}

\begin{proof}{}
We define a $\Ua$-module structure on 
$\Hom\left(V(\lambda),\Q(q)\right)$ by
\begin{equation*}
   \langle x f, u\rangle \defeq \langle f, \psi(x)u \rangle, \quad
   x\in \Ua, f\in\Hom\left(V(\lambda),\Q(q)\right),
   u\in V(\lambda).
\end{equation*}
This defines a $\Ua$-module structure since $\psi\colon\Ua\to
\Ua^{\operatorname{opp}}$ is an algebra homomorphism.
Let $u^\lambda$ be the unique linear form such that
\begin{equation*}
   \langle u^\lambda, G(b)\rangle =
   \begin{cases}
      1 & \text{if $G(b) = u_\lambda$}, \\
      0 & \text{otherwise}.
   \end{cases}
\end{equation*}
Then $u^\lambda$ has a weight $\lambda$. We claim that $u^\lambda$ is
an extremal vector.
From the definition all elements in a weight space
$\Hom\left(V(\lambda),\Q(q)\right)_\xi$ vanish on $V(\lambda)_\xi$.
Since weights of $V(\lambda)$ are contained in the convex hull of
$W\lambda$ \cite[Theorem~5.3]{Kas00}, the weights of $V'(\lambda)$
also have the same property. Therefore $u^\lambda$ is an extremal
vector. Now we have a $\Ua$-algebra homomorphism \( V(\lambda)\to
V'(\lambda)\subset\Hom\left(V(\lambda),\Q(q)\right) \) sending
$u_\lambda$ to $u^\lambda$. This defines a bilinear form satisfying
the desired properties. The uniqueness follows from the uniqueness of
the above homomorphism.
\end{proof}

\begin{Remark}\label{rem:bil2'}
The uniqueness holds even if \eqref{eq:bil2} holds only for
$x\in\Uq'(\ag)$.
In fact, this condition together with \eqref{eq:bil1} automatically
implies \eqref{eq:bil2} for $x = q^d$ as follows. When $u =
u_\lambda$, \eqref{eq:bil1} implies \eqref{eq:bil2} for $x =
q^d$. For a general case, we write $u = x u_\lambda$ with $x\in
\Uq'(\ag)_\xi$. Then
\begin{equation*}
\begin{split}
   & (q^d u, v) = q^{\langle d, \xi\rangle} (x q^d u_\lambda, v)
   = q^{\langle d, \xi\rangle}(q^d u_\lambda, \psi(x)v)
   = q^{\langle d, \xi\rangle}(u_\lambda, q^d \psi(x)v)
\\
   =\; & (u_\lambda, \psi(x) q^d v)
   = (x u_\lambda, q^d v)
   = (u, q^d v),
\end{split}
\end{equation*}
where we have used $\psi(x)\in \Uq'(\ag)_{-\xi}$.
\end{Remark}

\begin{Lemma}\label{lem:WeylAdj}
Let $M$ be an integrable $\Uq'(\ag)$-module with a bilinear form $(\
,\ )$ satisfying \eqref{eq:bil2} for $x\in \Uq'(\ag)$. Then
\begin{equation*}
   (T_w u, v) = (-1)^{N^\vee} q^{N} (u, T_{w^{-1}}v)
   \quad\text{for all $w\in\aW$, $u\in M_\xi$, $v\in M$},
\end{equation*}
where $N$ and $N^\vee$ are as in \lemref{lem:Weyl}.
\end{Lemma}

\begin{proof}{}
Let $T_{i,1}'$ be the operator defined in \cite[5.2.1]{Lu-Book}. A
direct calculation shows
\(
   (T_i u, v) = (u, T_{i,1}'v)
\)
for $u\in M_\xi$, $v\in M$. (We may assume that $v$ is contained in a
weight space. Thanks to \eqref{eq:bil2} for $x\in\Uq'(\ag)$, both hand
sides are $0$ unless the weight of $v$ is $s_i\xi + m\delta$ for some
$m\in\Z$.)
By [loc.\ cit., 5.2.3], we have
\(
  T_{i,1}' v = (-1)^{\langle h_i,\xi\rangle} q^{(\alpha_i,\xi)} T_i v
\).
The rest of the proof is the same as that of \lemref{lem:Weyl}.
\end{proof}

\begin{Lemma}\label{lem:WeylAdj2}
Let $M$ and $(\ ,\ )$ be as above. Let $u$, $v\in M$ be extremal
vectors. Then
\begin{equation*}
    (S_w u, v) = (u, S_{w^{-1}} v).
\end{equation*}
\end{Lemma}

\begin{proof}{}
Let $\xi$ be the weight of $u$. Using Lemmas~\ref{lem:WeylExtremal},
\ref{lem:WeylAdj}, we have
\begin{equation*}
   (S_w u, v)
   = (-1)^{N^\vee_+ + N^{\vee\prime}_+ + N^\vee} q^{-N_+ - N'_+ + N}
   (u, S_{w^{-1}}v),
\end{equation*}
where
\begin{equation*}
\begin{gathered}
     N = \sum_{\alpha\in\widehat\aR_+\cap w^{-1}(\widehat\aR_-)}
     (\alpha,\xi),
\qquad
     N_+ = \sum_{\alpha\in\widehat\aR_+\cap w^{-1}(\widehat\aR_-)}
     \max((\alpha,\xi),0),
\\
     N_+' = \sum_{\alpha'\in\widehat\aR_+\cap w(\widehat\aR_-)}
     \max((\alpha',w\xi),0),
\end{gathered}
\end{equation*}
and $N^\vee$, $N^\vee_+$, $N^{\vee\prime}_+$ are defined in similar
ways. Noticing $\alpha'\in\widehat\aR_+\cap
w(\widehat\aR_-)\Leftrightarrow - w^{-1}\alpha'\in \widehat\aR_+\cap
w^{-1}(\widehat\aR_-)$, we have $N = N_+ + N_+'$. Similarly we have
$N^\vee = N^\vee_+ + N^{\vee\prime}_+$. Therefore we have the
assertion.
\end{proof}

In order to study $(\ ,\ )$ on $V(\lambda)$ we need to relate it to
a bilinear form on the tensor product module $\widetilde V(\lambda)$.

\begin{Lemma}\label{lem:bilinear}
We have $(z_i u, z_i v) = (u,v)$ for $u,v\in V(\varpi_i)$.
\end{Lemma}

\begin{proof}{}
By the uniqueness, it is enough to show that $(z_i u, z_i v)$
satisfies (\ref{eq:bil1}, \ref{eq:bil2}). The property \eqref{eq:bil2} 
is clear. If $x\in \Uq'(\ag)$, then it holds since $z_i$ is
$\Uq'(\ag)$-linear. It also holds for $x = q^d$ thanks to $z_i q^d
z_i^{-1} = q^{-a_0d_i}q^d$.

Let us check \eqref{eq:bil1}. Since $\dim V(\varpi_i)_{\varpi_i} = 1$
by \cite[Proposition~5.10]{Kas00}, it is enough to show that
\(
   (z_i u_{\varpi_i}, z_i u_{\varpi_i}) = 1
\).
But this follows from the previous lemma.
\end{proof}

We define a $\Q(q)[z_i^\pm]$-valued bilinear form $(\!(\ ,\ )\!)$ on
$V(\varpi_i)$ by
\begin{equation*}
   (\!(u, v)\!) = 
   \begin{cases}
     z_i^{m} (z_i^{-m} u, v)
     & \text{if $\wt(u) = \wt(v) + md_i\delta$ for $m\in\Z$},
     \\
     0 & \text{otherwise}.
   \end{cases}
\end{equation*}
Since $z_i$ is $\Uq'(\ag)$-linear, we have
\begin{equation*}
   (\!(x u, v)\!) = (\!(u, \psi(x) v)\!)
   \quad\text{for $x\in\Uq'(\ag)$, $u,v\in V(\varpi_i)$.}
\end{equation*}
By \lemref{lem:bilinear} we have
\begin{equation}\label{eq:lsl}
   (\!( z_i^m u_{\varpi_i}, z_i^n u_{\varpi_i})\!)
   = z_i^{m-n}.
\end{equation}

We define a $\Q(q)[z_{i,\nu}^\pm]_{i\in I, \nu=1,\dots,m_i}$-valued
bilinear form $(\!(\ ,\ )\!)$ on $\widetilde V(\lambda)$ by
\begin{equation*}
   (\!(u, v)\!)
   \defeq \prod_{i,\nu} (\!(u_{i,\nu}, v_{i,\nu} )\!),
\end{equation*}
where $u_{i,\nu}$, $v_{i,\nu}$  is the $\nu$-th $V(\varpi_i)$-factor
of $u,v\in \widetilde V(\lambda)$.
We define a bilinear form $(\ ,\ )^\sim$ on $\widetilde V(\lambda)$ by
\begin{equation*}
   (u, v)^\sim \defeq \prod_{i\in I} \frac1{m_i!}
   \left[(\!(u, v)\!) \prod_{\mu\neq \nu} (1 - z_{i,\mu}z_{i,\nu}^{-1})
   \right]_1,
\end{equation*}
where $[f]_1$ denote the constant term in $f$.

\begin{Lemma}\label{lem:Schur}
Let $\bc_0$, $\bc'_0\in(\Z_{\ge 0}^{\aR_0})(\lambda)$. Then
\(
  (s_{\bc_0}(z)\widetilde u_\lambda,
    s_{\bc'_0}(z) \widetilde u_\lambda)^\sim = \delta_{\bc_0,\bc'_0}
\).
\end{Lemma}

\begin{proof}{}
Let $f = f(z)$ and $g = g(z)$ be polynomials in $z_{i,\nu}$'s ($i\in
I$, $\nu=1,\dots, m_i$). By \eqref{eq:lsl} we have
\begin{equation*}
   (f(z) \widetilde u_\lambda, g(z) \widetilde u_\lambda)^\sim
   = \prod_{i\in I} \frac1{m_i!}
   \left[f\overline{g}\prod_{\mu\neq \nu} (1 - z_{i,\mu}z_{i,\nu}^{-1})
   \right]_1,
\end{equation*}
where $\overline{g} = g(\dots, z_{i,\nu}^{-1},\dots)$. Considered as a 
bilinear form on the Laurent polynomial ring, it coincides with one in 
\cite[Chap.VI, \S9]{Mac} with $q=t$. The Schur functions give an
orthogonal base with respect to that bilinear form. Therefore we have
the assertion.
\end{proof}

\begin{Proposition}\label{prop:bilin}
Let $u,v\in V(\lambda)$. Then
\(
    (u,v) = (\Phi_\lambda(u), \Phi_\lambda(v))^\sim
\).
\end{Proposition}

\begin{proof}{}
It is enough to show that $(\Phi_\lambda(u), \Phi_\lambda(v))^\sim$
satisfies conditions in \propref{prop:bilinear}.
It is clear that the condition \eqref{eq:bil2} holds for
$x\in\Uq'(\ag)$. By \remref{rem:bil2'}, it is enough to check
\eqref{eq:bil1}. From \eqref{eq:bil2} for $x\in\Uq'(\ag)$, it is
enough to check \eqref{eq:bil1} when
$\operatorname{cl}(\operatorname{wt}(b)) =
\operatorname{cl}(\lambda)$, i.e., $\operatorname{wt}(b) = \lambda +
m\delta$ for some $m\in \Z$. Since weights of $V(\lambda)$ is
contained in the convex hull of $W\lambda$, $b$ is an extremal
vector. We have
\begin{equation*}
    (\Phi_\lambda(u_\lambda), \Phi_\lambda(G(b)))^\sim 
    = (\Phi_\lambda(S_w u_\lambda), \Phi_\lambda(S_w G(b)))^\sim
\end{equation*}
by \lemref{lem:WeylAdj2}.
We take $S_w$ as sufficiently many compositions of $\tf_i^{\max}$, we
may assume $S_w u_\lambda = S_{\bc_0}^- u_\lambda$, $S_w G(b) =
S_{\bc'_0}^-u_\lambda$.
(Recall that $S_{\bc_0}^- u_\lambda$ is an element of the global basis 
as we explained in \remref{rem:exgl}.)
Then
\begin{equation*}
    \left(\Phi_\lambda(u_\lambda), \Phi_\lambda(G(b))\right)^\sim 
    = (s_{\bc_0} \widetilde u_\lambda,
                 s_{\bc'_0} \widetilde u_\lambda)^\sim
    = \delta_{\bc_0, \bc_0'}
    = \delta_{u_\lambda, G(b)},
\end{equation*}
where we have used \propref{prop:tensor} and \lemref{lem:Schur}.
\end{proof}

From the proof of \propref{prop:bilin} the bilinear form $(\ ,\ )$ on
$V(\lambda)$ defined in \propref{prop:bilinear} also has the following 
characterization: it satisfies \eqref{eq:bil2} and
\(
    (S_{\bc_0} u_\lambda, S_{\bc'_0} u_\lambda) = \delta_{\bc_0,\bc_0'}.
\)
Since these conditions are symmetric, we have the following:

\begin{Corollary}\label{cor:symmetric}
The bilinear form $(\ ,\ )$ on $V(\lambda)$ is symmetric, i.e.,
$(u, v) = (v,u)$.
\end{Corollary}


\begin{Proposition}\label{prop:orth}
\textup{(1)} $(\La(\lambda),\La(\lambda))\subset \A_0$.

\noindent
Let $(\ ,\ )_0$ be the $\Q$-valued bilinear form on
$\La(\lambda)/q\La(\lambda)$ induced by $(\ ,\ )|_{q=0}$ on
$\La(\lambda)$.

\textup{(2)} $(\te_i u, v)_0 = (u, \tf_i v)_0$ for
$u,v\in \La(\lambda)/q\La(\lambda)$.

\textup{(3)} $\B(\lambda)$ is an orthonormal base with respect to $(\ ,
\ )_0$. In particular, $(\ ,\ )_0$ is positive definite.

\textup{(4)} 
\(
   \La(\lambda) = \{ u\in V\mid (u,u)\in \A_0 \}
\).
\end{Proposition}

\begin{proof}{}
We shall prove
\begin{itemize}
\item there exist representatives $\widetilde b$ for all
$b\in\B(\lambda)_\xi\subset \La(\lambda)_\xi/q\La(\lambda)_\xi$
such that
\(
   (\widetilde b, \widetilde b') \equiv \delta_{bb'} 
   \mod q\A_0
\)
for $b,b'\in \B(\lambda)_\xi$
\end{itemize}
by the induction on $(\xi, \xi)$.
Since $\La(\lambda)_\xi$ is spanned by $\widetilde b$'s over $\A_0$,
this implies the above equations for {\it any\/} representatives
$\widetilde b$. It also implies (1) and (3).
Recall
\(
   (\te_i\widetilde b, \widetilde b')
   = (\widetilde b, \tf_i\, \widetilde b')
\)
by \eqref{eq:psi_te}. Therefore the above assertion also implies (2).

First suppose that $b$ is extremal. Since we may assume that $\wt(b) =
\wt(b')$ by \eqref{eq:bil2}, we may assume $b'$ is also extremal by
\cite[5.3]{Kas00}. Then we may assume $\widetilde b =
S_{\bc_0}u_\lambda$, $\widetilde b' = S_{\bc'_0}u_\lambda$ by applying
$S_w$ for some $w\in \widehat W$. But, in this case, the assertion has
been already shown in \lemref{lem:Schur} and \propref{prop:bilin}.

Now we start the induction. Recall that $(\xi,\xi)$ is bounded from
above and $b\in\B(\lambda)$ is extremal if $(\wt b,\wt b)$ is maximal
(\cite[\S9.3]{Kas94}). Therefore when $(\xi,\xi)$ is maximal, both $b$
and $b'$ are extremal. We have already proved the assertion this case.

Now assuming the above for $\xi$ such that $(\xi,\xi) > a$, let us
prove it for $\xi$ with $(\xi,\xi) = a$. For $i\in I$, suppose that
$\langle h_i,\xi\rangle \ge 0$. We consider $\te_i b$. If $\te_i b\neq
0$, then we have
\begin{equation*}
   \left(\wt(\te_i b),\wt(\te_i b)\right) = (\xi+\alpha_i,\xi+\alpha_i)
   > (\xi,\xi).
\end{equation*}
Therefore we have
\begin{equation*}
   \left(\tf_i \te_i \widetilde b, \widetilde b'\right)
   = \left( \te_i\widetilde b, \te_i\widetilde b'\right)
   \equiv \delta_{\te_i b,\te_i b'}
   \equiv \delta_{bb'} \mod q\A_0
\end{equation*}
by the induction hypotheis. Hence the assertion holds if we replace
the representative $\widetilde b$ by another representative $\tf_i
\te_i\widetilde b$.
Similarly, if $\langle h_i,\xi\rangle \le 0$ and $\tf_i b\neq 0$, we
replace $\widetilde b$ by $\te_i\tf_i \widetilde b$ to get the
assertion.

Since we may suppose that $b$ is not extremal, there exists $w\in\aW$
such that $S_w b$ satisfies $\te_i S_w b \neq 0$ if $\langle h_i,
w\xi\rangle \ge 0$ and $\tf_i S_w b\neq 0$ if $\langle h_i,
w\xi\rangle \le 0$. Then we have $(\tf_i\te_i S_w \widetilde b,
S_w\widetilde b')$ or $(\te_i\tf_i S_w\widetilde b, S_w\widetilde b')$
is in $\delta_{bb'} + q\Z[q]$. Therefore we are done.

The statement (4) follows from \cite[14.2.2]{Lu-Book}.
\end{proof}

The followign result generalizes \cite[Theorem A]{VV2} from fundamental
representations to arbitrary $\lambda$:

\begin{MainTheorem}\label{thm:main2}
\textup{(1)} $\{ G(b) \}_{b\in \B(\lambda)}$ is almost orthonormal for 
$(\ ,\ )$, that is, $(G(b), G(b')) \equiv \delta_{bb'} \mod
q\Z[q]$.

\textup{(2)}
\(
   \{ \pm G(b) \mid b\in\B(\lambda)\}
   = \left\{ u\in V^\Z(\lambda) \left|\, \overline{u} = u, \;
       (u,u) \equiv 1 \bmod q\Z[q] \right\}\right..
\)
\end{MainTheorem}

\begin{proof}
We claim
\begin{equation*}
   (u, v)\in \Z[q,q^{-1}] \quad\text{for $u,v\in V^\Z(\lambda)$}.
\end{equation*}
The assertion is obvious for the special case $u = u_\lambda$ by
\eqref{eq:bil1}. For general case, we may assume $u = x u_\lambda$ for
$x \in \Ui$. Then
\(
   (xu_\lambda, v) = (u_\lambda, \psi(x) v)
\).
Since $\psi(x)\in\Ui$ and $V^\Z(\lambda)$ is stable under the action
of $\Ui$, the assertion follows from the special case.

Combining with \propref{prop:orth}, we have
\begin{equation*}
   (G(b), G(b')) - \delta_{bb'} \in \Z[q,q^{-1}]\cap q\A_0
   = q\Z[q].
\end{equation*}
This is the statement~(1). The statement~(2) follows from
the argument of \cite[14.2.3]{Lu-Book}.
\end{proof}

\begin{Remark}
%
%
%
  Lusztig conjectures that the universal standard module $M(\lambda)$,
  more precisely its tensor product of\linebreak[1]
  $\otimes_{R(G_\lambda)} R(H_\lambda)$, which is isomorphic to
  $\Breve V^\Z(\lambda)$, has a signed base characterized by
  the almost orthogonality property \thmref{thm:main2}(2), with
  respect to geometrically defined bilinear form and bar involution
  \cite{Lu-rem}. (See \subsecref{subsec:standard} for notations.)
  Recently Varagnolo-Vasserot~\cite{VV2} give a proof of the
  conjecture by showing that $\{ G(b) \mid b\in \Breve\B(\lambda)\}$
  satisfies the property.
%
%
They also conjecture that the global base $\{ G(b) \mid
b\in\B(\lambda)\}$ of $V(\lambda)$ satisfies the almost orthogonality
property with respect to the geometric bilinear form and bar
involution.
Their conjecture follows from \thmref{thm:main2}(2) since the
geometric bilinear form and bar involution coincide with ones used in
this paper, as Varangnolo and Vasserot proved that the formers
satisfy the conditions in \propref{prop:bilinear} (more precisely
\eqref{eq:bil2} and \( (S_{\bc_0} u_\lambda, S_{\bc'_0} u_\lambda) =
\delta_{\bc_0,\bc_0'} \)) and the equality $\overline{xu} =
\overline{x}\ \overline{u}$. Remark that these hold only after an {\it
appropriate\/} normalization of standard modules so that we have
$x_{i,\nu} = \pm z_{i,\nu}$. This is the normalization in \cite{VV2}
different from ours. This point is clarified during discussion with
Varagnolo-Vasserot in Februrary 2002.
\end{Remark}

\end{document}